# THE CONTACT STRUCTURE ON THE LINK OF A CUSP SINGULARITY

NAOHIKO KASUYA

ABSTRACT. In this paper, we study the contact structures on the links of cusp singularities. We show that they are contactomorphic to Sol-manifolds with the positive contact structures arising from Anosov flows. As an application, we obtain an exotic contact structure on $S^5$ by performing Lutz-Mori twist along the Sol-manifold.

## 1. INTRODUCTION

The algebraic surface $V$ in $\mathbb{C}^3$ defined by the equation $f(x,y,z) = 0$, where

$$f(x,y,z) = x^p + y^q + z^r + xyz \text{ with } p,q,r \in \mathbb{Z}_{\geq 2} \text{ satisfying } \frac{1}{p} + \frac{1}{q} + \frac{1}{r} < 1,$$

has the only singularity at $(0,0,0)$. This singularity is called a cusp singularity or a hyperbolic singularity. The canonical contact structure on the link $K = V \cap S_\varepsilon^5$, where $S_\varepsilon^5$ is a small sphere centered at $(0,0,0)$, is the restriction of the standard contact structure $(S^5, \xi_0)$. The link $K$ is diffeomorphic to a Sol-manifold by the result of Laufer[6]. On a Sol-manifold $T_A$, there is an Anosov flow and it induces the positive contact structure $(T_A, \beta_+ + \beta_-)$ and the negative contact structure $(T_A, \beta_+ - \beta_-)$, where $\beta_+$ and $\beta_-$ are the 1-forms defining the Anosov foliations. It is the bi-contact structure on $T_A$. Our main theorem is following.

**Theorem 1.1.** *The link $(K, \xi_0|K)$ is contactomorphic to the Sol-manifold $(T_A, \beta_+ + \beta_-)$, where*

$$A = \left\{ \begin{pmatrix} 0 & 1 \\ -1 & p \end{pmatrix} \begin{pmatrix} 0 & 1 \\ -1 & q \end{pmatrix} \begin{pmatrix} 0 & 1 \\ -1 & r \end{pmatrix} \right\}^{-1}.$$

The strategy for the proof is following. The cusp singularity $(V, (0,0,0))$ is analytically equivalent to the cusp singularity $\infty$ in the analytic surface $\overline{\mathbb{H} \times \mathbb{H}/G(M,V)}$ which Hirzebruch constructed in [3]. Therefore it is enough to prove that the canonical contact structure on the link of $\infty$ is contactomorphic to the Sol-manifold $(T_A, \beta_+ + \beta_-)$. We constructed the contactomorphism explicitly in Theorem 3.15. One can also prove this result by Honda's classification of tight contact structures on $T_A$ [4] and the theorem of Y.Lekili and B.Ozbagci that Milnor fillable contact structures are universally tight [7] (see Remark 3.19).

Atsuhide Mori also constructed examples of surface singularities in $\mathbb{C}^3$ whose links are contactomorphic to Sol-manifolds in theorem 3.5 of [11]. He noticed the supporting open book decomposition of $(T_A, \beta_+ + \beta_-)$ and realized it as the Milnor open book of the link of a surface singularity in $\mathbb{C}^3$. His examples are corresponding to the cases of theorem 1.1 where $p = 2$. His result is in the context of constructing the examples of Lutz-Mori twists, the generalization of Lutz twists to higher dimensions. We can modify the standard contact sphere $(S^5_\varepsilon, \xi_0)$ by the Lutz-Mori twist along the link $K$. The resultant contact structure $(S^5_\varepsilon, \xi)$ is equivalent to the contact structure $\ker(\alpha_{m,k})$ ($m = 1, 2, 3$) of Theorem 3.7 in [11], which is an exotic contact structure on $S^5$ and can be deformed via contact structures to a spinnable foliation (see the section 4 and [11]).

The author thanks Professor Takashi Tsuboi, Atsuhide Mori, Hiroki Kodama, Ko Honda, Patrick Massot and Burak Ozbagci for many helpful advices.

## 2. Sol-manifolds and bi-contact structures

Let $\begin{pmatrix} x \\ y \end{pmatrix}$ be the coordinate on the torus $T^2 = \mathbb{R}^2/\mathbb{Z}^2$ and $(\begin{pmatrix} x \\ y \end{pmatrix}, z)$ be the coordinate on $T^2 \times [0, 1]$. Let $A$ be an element of $SL(2, \mathbb{Z})$ such that $\operatorname{tr}(A) > 2$. Then $A$ has two positive eigenvalues $a$ and $a^{-1}$ and the corresponding eigenvectors $v_+$ and $v_-$, where $a > 1$ and $dx \wedge dy(v_+, v_-) = 1$.

**Definition 2.1.** *We define an equivalence relation $\sim$ on $T^2 \times [0, 1]$ by $(A \begin{pmatrix} x \\ y \end{pmatrix}, 0) \sim (\begin{pmatrix} x \\ y \end{pmatrix}, 1)$. The quotient $T_A = T^2 \times [0, 1]/\sim$ is called a hyperbolic mapping torus.*

**Definition 2.2.** *The Lie group $Sol^3$ is the split extension $1 \to \mathbb{R}^2 \to Sol^3 \to \mathbb{R} \to 1$ whose group structure is given by*

$$(u, v; w) \cdot (u', v'; w') = (u + e^w u', v + e^{-w} v'; w + w') \text{ on } \mathbb{R}^2 \times \mathbb{R}.$$

*There is a left invariant metric $e^{-2w} du \otimes du + e^{2w} dv \otimes dv + dw \otimes dw$ on $Sol^3$. Let $\Gamma$ be a cocompact discrete subgroup of $Sol^3$. The compact quotient $M^3 = \Gamma \backslash Sol^3$ is called a Sol-manifold.*

The discrete subgroups of $Sol^3$ are all of the forms $\Gamma = M \rtimes V = (\mathbb{Z} \times \mathbb{Z}) \rtimes \mathbb{Z}$. The quotient of $\mathbb{R}^2$ by the lattice $M$ is $T^2$, while the quotient of $\mathbb{R}$ by $V$ is a circle. The quotient $M^3 = \Gamma \backslash Sol^3$ is a $T^2$-bundle over $S^1$ with hyperbolic monodromy. Conversely, a hyperbolic mapping torus $T_A$ is a Sol-manifold. From now on, $T_A$ represents a Sol-manifold. The left invariant 1-forms $e^{-w} du$ and $-e^w dv$ on $Sol^3$ induce the 1-forms $\beta_+ = a^{-z} dx \wedge dy(v_+, \cdot)$ and $\beta_- = -a^z dx \wedge dy(v_-, \cdot)$ on $T_A$.

**Definition 2.3.** *A non-singular flow $\phi_t$ on $M^3$ is an Anosov flow if the tangent bundle $TM^3$ has the $\phi_t$-invariant decomposition for some Riemannian metric $g$ on $M^3$ such that $TM^3 = T\phi \oplus E^{uu} \oplus E^{ss}$, where*

$$T\phi = \{\text{tangent vectors along the flow lines}\},$$

$$E^{uu} = \left\{v \in TM^3; ||T\phi_t(v)|| \geq \exp(ct)||v||, t > 0\right\},$$

$$E^{ss} = \left\{v \in TM^3; ||T\phi_t(v)|| \geq \exp(ct)||v||, t < 0\right\},$$

*for some positive real number $c$. We call $E^s = T\phi \oplus E^{ss}$ (resp. $E^u = T\phi \oplus E^{uu}$) weakly stable (resp. unstable) plane field. We obtain two codimension 1 foliations, the unstable foliation $F^u$ and the stable foliation $F^s$ as the tangent bundles of $E^u$ and $E^s$ respectively. They are called Anosov foliations.*

**Definition 2.4.** *A bi-contact structure on $M^3$ is a pair of a positive contact structure $\xi_+$ and a negative contact structure $\xi_-$ on $M^3$ which are transverse to each other.*

**Example 2.5.** *Let $M^3$ be a Sol-manifold $T_A$ and fix a Riemannian metric $g = \beta_+ \otimes \beta_+ + \beta_- \otimes \beta_- + dz \otimes dz$. Then the flow $\phi_t(x, y, z) = (x, y, z+t)$ is an Anosov flow with the stable foliation $F^s = \ker \beta_+ = <(\frac{\partial}{\partial z}), v_->$ and the unstable foliation $F^u = \ker \beta_- = <(\frac{\partial}{\partial z}), v_+>$. Moreover, $\xi_+ = \ker(\beta_+ + \beta_-)$ and $\xi_- = \ker(\beta_+ - \beta_-)$ form a bi-contact structure on $T_A$. We can see that the flow $\phi_t$ pushes the plane fields $\xi_+$ and $\xi_-$ towards $F^u$ (or $F^s$ if you flow backward). More precisely,*

$$\lim_{t \to +\infty} (\phi_t)_* \xi_+ = \lim_{t \to +\infty} (\phi_t)_* \xi_- = E^u$$

*and*

$$\lim_{t \to -\infty} (\phi_t)_* \xi_+ = \lim_{t \to -\infty} (\phi_t)_* \xi_- = E^s.$$

See [10] for more about Anosov foliations and bi-contact structures.

**Remark 2.6.** *The 1-forms $\beta_+ + \beta_-$ and $\beta_+ - \beta_-$ are induced by left invariant contact forms $e^{-w}du - e^w dv$ and $e^{-w}du + e^w dv$ on $Sol^3$, respectively. The universal covering of $(T_A, \beta_+ + \beta_-)$ is $(Sol^3, e^{-w}du - e^w dv)$, which is the standard positive contact structure on $\mathbb{R}^3$. Therefore $(T_A, \beta_+ + \beta_-)$ is universally tight and it is the unique positive contact structure on $T_A$ that is universally tight and minimally twisting [4]. Similarly, $(T_A, \beta_+ - \beta_-)$ is the unique negative contact structure on $T_A$ that is universally tight and minimally twisting.*

## 3. Cusp singularities and the main theorem

The link of a cusp singularity constructed by the method of Hirzebruch[3] is diffeomorphic to a Sol-manifold. Moreover, this link has a canonical contact structure which is equivalent to $(T_A, \beta_+ + \beta_-)$. By Laufer and Karras's criterion on the embeddability of a cusp singularity, this link can be embedded into $\mathbb{C}^3$ as the link of the surface singularity $T_{p,q,r}$ for some $p, q, r \in \mathbb{Z}_{\geq 2}$. The main theorem is that the restriction of the standard contact srtucture $(S_\varepsilon^5, \xi_0)$ to the embedded link of $T_{p,q,r}$ is contactomorphic to $(T_A, \beta_+ + \beta_-)$.

### 3.1. Cusp singularities.

**Definition 3.1.** *Let $(V, 0)$ be an isolated surface singularity germ. $(V, 0)$ is normal if every bounded holomorphic function on $V^* = V - 0$ extends to a holomorphic function at $0$.*

**Definition 3.2.** *Let $(V, 0)$ be a normal surface singularity. Then there exists a non-singular complex surface $M$ and a proper analytic map $\pi : M \to V$ satisfying the following conditions (1) and (2).*
  (1) *$E = \pi^{-1}(0)$ is a divisor in $M$, i.e., a union of 1-dimensional compact curves in $M$; and*
  (2) *The restriction of $\pi$ to $\pi^{-1}(V^*)$ is a biholomorphic map between $M - E$ and $V^*$.*

*The surface $M$ is called a resolution of the singularity of $V$, and $\pi : M \to V$ is the resolution map. The divisor $E$ is called the exceptional set. The divisor $E$ is good if it satisfies the following two conditions:*

  (3) *each irreducible component $E_i$ of $E$ is non-singular; and*
  (4) *$E$ has normal crossings, i.e., $E_i$ intersects $E_j$, $i \neq j$, in at most one point, where they meet transversally, and no three of them intersect.*

**Definition 3.3.** *A resolution $\pi : M \to V$ is minimal if any resolution $\pi' : M' \to V$, there is a proper analytic map $p : M' \to M$ such that $\pi' = \pi \circ p$.*

By Castelnuovo's criterion, minimality of a resolution is equivalent to the condition that the exceptional set contains no non-singular rational curves with self-intersection $-1$.

**Definition 3.4.** *Let $(V, 0)$ be a normal surface singularity and $\pi : M \to V$ be the minimal resolution of $(V, 0)$. $(V, 0)$ is called a cusp singularity if the exceptional set $E = \pi^{-1}(0)$ is a cycle of non-singular rational curves, i.e., $E_i$ intersects $E_{i+1}$ transversally at one point for $1 \leq i \leq n$ and $E$ has no other crossings, where $E = \bigcup_{i=1}^{n} E_i$ and $E_{n+1}$ means $E_1$.*

If $(V, 0)$ is a cusp singularity embedded in $\mathbb{C}^3$, the link $K = V \cap S_\varepsilon^5$ is a Sol-manifold. The following theorem gives embeddings of cusp singularities into $\mathbb{C}^3$.

**Theorem 3.5** (Karras[5], Laufer[6]). *A cusp singularity is embedded in $\mathbb{C}^3$ if and only if it is defined by the equation $x^p + y^q + z^r + xyz = 0$ ($\frac{1}{p} + \frac{1}{q} + \frac{1}{r} < 1$). The singularity $(0, 0, 0)$ is called $T_{p,q,r}$ singularity. The link $K$ of $T_{p,q,r}$ is diffeomorphic to the Sol-manifold $T_A$, where*

$$A = \left\{ \begin{pmatrix} 0 & 1 \\ -1 & p \end{pmatrix} \begin{pmatrix} 0 & 1 \\ -1 & q \end{pmatrix} \begin{pmatrix} 0 & 1 \\ -1 & r \end{pmatrix} \right\}^{-1}.$$

The canonical contact structure on the link $K$ is the restriction $(K, \xi_0|K)$ of the standard contact structure $(S_\varepsilon^5, \xi_0)$. To determine this contact structure, let us look at the construction of cusp singularities by Hirzebruch[3] in the next section.

3.2. **Hirzebruch's construction.** We are going to construct cusp singularities by the method of Hirzebruch[3] until Lemma 3.14. Let us assume that $\mathbb{Z} \ni k \mapsto b_k \in \mathbb{N}$ is a function satisfying $b_k \geq 2$ for all $k$, and $b_k \geq 3$ for some $k$. For each integer $k$, take a copy $R_k$ of $\mathbb{C}^2$ with coordinates $(u_k, v_k)$. We define $R'_k = R_k \setminus \{u_k = 0\}$ and $R''_k = R_k \setminus \{v_k = 0\}$. The equations

(1) $$u_{k+1} = u_k^{b_k} v_k,$$

(2) $$v_{k+1} = \frac{1}{u_k}$$

give a biholomorphic map $\varphi_k : R'_k \to R''_{k+1}$. In the disjoint union $\bigcup R_k$ we make all the identifications (1) and (2). We get a set $Y$. Denote the coordinate map $R_j \to \mathbb{C}^2$ by $\psi_j$. The map $\psi_k \circ \psi_0^{-1}$ is given by

$$u_k = u_0^{p_k} v_0^{q_k},$$
$$v_k = \frac{1}{u_0^{p_{k-1}} v_0^{q_{k-1}}},$$

where

$$\begin{pmatrix} p_k & q_k \\ -p_{k-1} & -q_{k-1} \end{pmatrix} = \begin{pmatrix} b_{k-1} & 1 \\ -1 & 0 \end{pmatrix} \begin{pmatrix} b_{k-2} & 1 \\ -1 & 0 \end{pmatrix} \cdots \begin{pmatrix} b_0 & 1 \\ -1 & 0 \end{pmatrix},$$

$$\frac{p_k}{q_k} = b_0 - \cfrac{1}{b_1 - \cfrac{1}{b_2 - \cdots - \cfrac{1}{b_{k-1}}}}.$$

Moreover $\{v_k = 0\}$ and $\{u_{k+1} = 0\}$ are pasted together by the equation (2) and form $\mathbb{C}P^1$, and the self-intersection number of $\mathbb{C}P^1$ is $-b_k$ by the equation (1). Let us put this rational curve $S_k$, then the intersection numbers are

$$S_i \cdot S_i = -b_i, \ S_i \cdot S_{i+1} = 1, \ S_i \cdot S_j = 0 \ (\mid i - j \mid \geq 2).$$

Let us assume that the function $\{b_k\}$ above is periodic and $r$ is the period: $b_{k+r} = b_k (k \in \mathbb{Z})$. We consider the infinite continued fractions

$$w_k = b_k - \cfrac{1}{b_{k+1} - \cfrac{1}{b_{k+2} - \cdots}}.$$

**Proposition 3.6.** $w_k$ *satisfies the following conditions.*

(1) $w_k$ *is a quadratic irrational number which is greater than* 1.
(2) $w_k = w_{k+r}$.

*Proof.* (1) If $b_j = 2$ for all $j$, $w_k$ are all equal to 1. Hence $w_k$ is greater than 1 by the assumption that $b_i \geq 2$ for all $i$ and $b_i \geq 3$ for some $i$. Now we will show that $w_k$ is a quadratic irrational number. It is enough to show that $w_0$ is so. Since

$$\frac{p_k}{q_k} = \frac{b_{k-1}p_{k-1} - p_{k-2}}{b_{k-1}q_{k-1} - q_{k-2}} = b_0 - \cfrac{1}{b_1 - \cfrac{1}{b_2 - \cdots - \cfrac{1}{b_{k-1}}}}$$

and

$$w_0 = b_0 - \cfrac{1}{b_1 - \cfrac{1}{b_2 - \cdots - \cfrac{1}{b_{r-1} - \cfrac{1}{w_0}}}},$$

we get

$$w_0 = \frac{(b_{r-1} - \frac{1}{w_0})p_{r-1} - p_{r-2}}{(b_{r-1} - \frac{1}{w_0})q_{r-1} - q_{r-2}} = \frac{p_r - \frac{1}{w_0}p_{r-1}}{q_r - \frac{1}{w_0}q_{r-1}} = \frac{w_0 p_r - p_{r-1}}{w_0 q_r - q_{r-1}}.$$

Hence $w_0$ is a solution of a quadratic equation $q_r t^2 - (p_r + q_{r-1})t + p_{r-1} = 0$. Moreover it is not a rational number since continued fraction of a rational number has a finite length or the repeating section [2]. Thus $w_0$ is a quadratic irrational number.

(2) It is clear by the condition $b_i = b_{i+r}$ for all $i$ and the definition of $w_k$. □

The quadratic irrational numbers $w_k$ all belong to the same real quadratic field $K = \mathbb{Q}[w_0]$. We consider the $\mathbb{Z}$-module $M = \mathbb{Z} \cdot w_0 + \mathbb{Z} \cdot 1$.

We define the action on $\mathbb{C}^2$ as follows. For $a \in M$, $\bar{a}$ denotes the conjugate irrational number of the quadratic irrational number $a$.

$$a : (z_1, z_2) \mapsto (z_1 + a, z_2 + \bar{a}).$$

Then $\mathbb{C}^2/M$ is diffeomorphic to $T^2 \times \mathbb{R}^2$. We show the map

$$\Phi : Y - \bigcup_{j \in \mathbb{Z}} S_j \to \mathbb{C}^2/M$$

given by $\Phi(u_0, v_0) = (z_1, z_2)$, where

$$2\pi i z_1 = w_0 \log u_0 + \log v_0, \quad 2\pi i z_2 = \bar{w}_0 \log u_0 + \log v_0$$

is well-defined and biholomorphic.

**Proposition 3.7.** *The map $\Phi$ is well-defined and biholomorphic.*

*Proof.* The logarithms are defined modulo $2\pi i$. For $k, l \in \mathbb{Z}$, we have

$$w_0(\log u_0 + 2\pi k i) + (\log v_0 + 2\pi l i) = 2\pi i(z_1 + (kw_0 + l)) \text{ and}$$
$$\bar{w}_0(\log u_0 + 2\pi k i) + (\log v_0 + 2\pi l i) = 2\pi i(z_2 + \overline{(kw_0 + l)}).$$

Since $kw_0 + l$ is in $M$, $(z_1, z_2)$ is defined modulo $M$. So the map $\Phi$ is well-defined. Next we construct $\Psi$, the inverse map of $\Phi$. We solve the equation

$$\begin{cases} 2\pi i z_1 = w_0 \log u_0 + \log v_0 \\ 2\pi i z_2 = \bar{w}_0 \log u_0 + \log v_0 \end{cases}$$

for $\log u_0$ and $\log v_0$ and we get

$$\begin{pmatrix} \log u_0 \\ \log v_0 \end{pmatrix} = \frac{2\pi i}{w_0 - \bar{w}_0} \begin{pmatrix} 1 & -1 \\ -\bar{w}_0 & w_0 \end{pmatrix} \begin{pmatrix} z_1 \\ z_2 \end{pmatrix},$$

where $(z_1, z_2)$ is given only modulo $M$. For any $a \in M$, we have $k, l \in \mathbb{Z}$ such that $a = kw_0 + l$, and

$$\begin{aligned} \frac{2\pi i}{w_0 - \bar{w}_0} \begin{pmatrix} 1 & -1 \\ -\bar{w}_0 & w_0 \end{pmatrix} \begin{pmatrix} z_1 + a \\ z_2 + \bar{a} \end{pmatrix} &= \frac{2\pi i}{w_0 - \bar{w}_0} \begin{pmatrix} 1 & -1 \\ -\bar{w}_0 & w_0 \end{pmatrix} \begin{pmatrix} z_1 + (kw_0 + l) \\ z_2 + (k\bar{w}_0 + l) \end{pmatrix} \\ &= \frac{2\pi i}{w_0 - \bar{w}_0} \begin{pmatrix} 1 & -1 \\ -\bar{w}_0 & w_0 \end{pmatrix} \begin{pmatrix} z_1 \\ z_2 \end{pmatrix} + 2\pi i \begin{pmatrix} k \\ l \end{pmatrix} \\ &= \begin{pmatrix} \log u_0 + 2\pi k i \\ \log v_0 + 2\pi l i \end{pmatrix}. \end{aligned}$$

Hence $\log u_0$ and $\log v_0$ are defined modulo $2\pi i$. Therefore $u_0$ and $v_0$ are uniquely defined. Thus the map

$$\Psi : \mathbb{C}^2/M \to Y - \bigcup_{j \in \mathbb{Z}} S_j$$

given by $\Psi(z_1, z_2) = (u_0, v_0)$, where

$$u_0 = \exp 2\pi i \left(\frac{z_1 - z_2}{w_0 - \bar{w}_0}\right), \quad v_0 = \exp 2\pi i \left(\frac{w_0 z_2 - \bar{w}_0 z_1}{w_0 - \bar{w}_0}\right),$$

is the inverse of $\Phi$. Since $\Psi$ is a composition of the logarithm function and linear functions, and $\Phi$ is a composition of the exponential function and linear functions, they are holomorphic functions. Therefore $\Phi$ is biholomorphic. □

Put $A_0 = 1$ and $A_{k+1} = w_{k+1}^{-1} A_k$. Then

$$A_k = (w_1 w_2 \cdots w_k)^{-1} (k \geq 1), \quad A_0 = 1, \quad A_k = w_{k+1} \cdots w_{-1} w_0 (k \leq -1).$$

**Proposition 3.8.** $A_k$ satisfies the following conditions.
  (1) $0 < A_{k+1} < A_k$.
  (2) $A_{k+1} = b_k A_k - A_{k-1}$.
  (3) $2\pi i z_1 = A_{k-1} \log u_k + A_k \log v_k$,
      $2\pi i z_2 = \bar{A}_{k-1} \log u_k + \bar{A}_k \log v_k$.
  (4) $A_{k+r} = A_r A_k$,
      $(A_r)^n = A_{nr}$.

*Proof.* (1) It is clear by $A_0 = 1$, $A_{k+1} = w_{k+1}^{-1} A_k$ and $w_k > 1$.
  (2) We get $A_{k-1} = b_k A_k - A_{k+1}$ by multiplying $A_k$ to the both sides of $w_k = b_k - \dfrac{1}{w_{k+1}}$. Hence $A_{k+1} = b_k A_k - A_{k-1}$.
  (3) Since $\log u_k = p_k \log u_0 + q_k \log v_0$ and $\log v_k = -p_{k-1} \log u_0 - q_{k-1} \log v_0$, it is enough to show that $A_{k-1} p_k - A_k p_{k-1} = w_0$ and $A_{k-1} q_k - A_k q_{k-1} = 1$.

$$\begin{pmatrix} p_k & -p_{k-1} \\ q_k & -q_{k-1} \end{pmatrix} \begin{pmatrix} A_{k-1} \\ A_k \end{pmatrix} = \begin{pmatrix} b_0 & -1 \\ 1 & 0 \end{pmatrix} \cdots \begin{pmatrix} b_{k-2} & -1 \\ 1 & 0 \end{pmatrix} \begin{pmatrix} b_{k-1} & -1 \\ 1 & 0 \end{pmatrix} \begin{pmatrix} A_{k-1} \\ A_k \end{pmatrix}$$

$$= \begin{pmatrix} b_0 & -1 \\ 1 & 0 \end{pmatrix} \cdots \begin{pmatrix} b_{k-2} & -1 \\ 1 & 0 \end{pmatrix} \begin{pmatrix} A_{k-2} \\ A_{k-1} \end{pmatrix}$$

$$= \begin{pmatrix} A_{-1} \\ A_0 \end{pmatrix}$$

$$= \begin{pmatrix} w_0 \\ 1 \end{pmatrix}.$$

(4) It is clear by $w_k = w_{k+r}$ and the definition of $A_k$. $\square$

**Proposition 3.9.** $A_r$ is a unit in $M$ which satisfies $A_r < 1$.

*Proof.* Since $\begin{pmatrix} A_{k+1} \\ A_k \end{pmatrix} = \begin{pmatrix} b_k & -1 \\ 1 & 0 \end{pmatrix} \begin{pmatrix} A_k \\ A_{k-1} \end{pmatrix}$,

$$\begin{pmatrix} A_{k+1} \\ A_k \end{pmatrix} = \begin{pmatrix} b_k & -1 \\ 1 & 0 \end{pmatrix} \cdots \begin{pmatrix} b_0 & -1 \\ 1 & 0 \end{pmatrix} \begin{pmatrix} A_0 \\ A_{-1} \end{pmatrix}.$$

Since $A_{k+r} = A_r A_k$,

$$A_r \begin{pmatrix} A_0 \\ A_{-1} \end{pmatrix} = \begin{pmatrix} A_r \\ A_{r-1} \end{pmatrix} = \begin{pmatrix} b_{r-1} & -1 \\ 1 & 0 \end{pmatrix} \cdots \begin{pmatrix} b_0 & -1 \\ 1 & 0 \end{pmatrix} \begin{pmatrix} A_0 \\ A_{-1} \end{pmatrix},$$

that is,

$$A_r \begin{pmatrix} 1 \\ w_0 \end{pmatrix} = \begin{pmatrix} b_{r-1} & -1 \\ 1 & 0 \end{pmatrix} \cdots \begin{pmatrix} b_0 & -1 \\ 1 & 0 \end{pmatrix} \begin{pmatrix} 1 \\ w_0 \end{pmatrix}.$$

Take the conjugate of both sides, we have

$$\bar{A}_r \begin{pmatrix} 1 \\ \bar{w}_0 \end{pmatrix} = \begin{pmatrix} b_{r-1} & -1 \\ 1 & 0 \end{pmatrix} \cdots \begin{pmatrix} b_0 & -1 \\ 1 & 0 \end{pmatrix} \begin{pmatrix} 1 \\ \bar{w}_0 \end{pmatrix}.$$

These show that the eigenvalues of $\begin{pmatrix} b_{r-1} & -1 \\ 1 & 0 \end{pmatrix} \cdots \begin{pmatrix} b_0 & -1 \\ 1 & 0 \end{pmatrix} \in SL(2;\mathbb{Z})$ are $A_r$ and $\bar{A}_r = A_r^{-1}$, and the corresponding eigenvectors are $\begin{pmatrix} 1 \\ w_0 \end{pmatrix}$ and $\begin{pmatrix} 1 \\ \bar{w}_0 \end{pmatrix}$. Since $\{1, w_0\}$ is a basis of $M$ and $\begin{pmatrix} b_{r-1} & -1 \\ 1 & 0 \end{pmatrix} \cdots \begin{pmatrix} b_0 & -1 \\ 1 & 0 \end{pmatrix} \in SL(2;\mathbb{Z})$, $\{A_r, A_r w_0\}$ is also a basis. Hence $A_r M = M$ and $A_r$ is a unit in $M$. Moreover it is clear that $A_r < 1$ by $A_r = (w_1 \cdots w_r)^{-1}$ and $w_k > 1$. □

We consider the group
$$V = \{(A_r)^n \mid n \in \mathbb{Z}\}.$$

We define the actions of $V$ as follows. Let $\mathbb{H}$ be the upper half plane, and put
$$Y^+ = \Phi^{-1}(\mathbb{H} \times \mathbb{H}/M) \cup \bigcup_{j \in \mathbb{Z}} S_j.$$

The action of $V$ on $Y^+$ is given by $(A_r)^n : R_k \to R_{k+nr}$, where
$$(A_r)^n(u_k, v_k) = (u_k, v_k) \text{ for } (A_r)^n \in V.$$

We define the action of $V$ on $\mathbb{H} \times \mathbb{H}/M$ by
$$(A_r)^n \cdot (z_1, z_2) = ((A_r)^n z_1, (\bar{A}_r)^n z_2) \text{ for } (A_r)^n \in V.$$

**Proposition 3.10.** *These two actions of $V$ on $Y^+$ and on $\mathbb{H} \times \mathbb{H}/M$ are equivariant with respect to $\Phi$, and the action of $V$ on $Y^+$ is free and properly discontinuous.*

*Proof.* By Proposition 4.7.(3)(4), we have
$$\begin{aligned}
2\pi i((A_r)^n z_1) &= (A_r)^n A_{k-1} \log u_k + (A_r)^n A_k \log v_k \\
&= A_{k+nr-1} \log u_k + A_{k+nr} \log v_k \text{ and} \\
2\pi i((\bar{A}_r)^n z_2) &= (\bar{A}_r)^n \bar{A}_{k-1} \log u_k + (\bar{A}_r)^n \bar{A}_k \log v_k \\
&= \bar{A}_{k+nr-1} \log u_k + \bar{A}_{k+nr} \log v_k.
\end{aligned}$$

Hence the coordinate of $\Psi((A_r)^n \cdot (z_1, z_2))$ in $R_{k+nr}$ is $(u_k, v_k)$. Therefore these two actions of $V$ on $Y^+$ and on $\mathbb{H} \times \mathbb{H}/M$ are equivariant with respect to $\Phi$. Next we show that the action of $V$ on $Y^+$ is free. First we show that the action is free on $Y^+ - \bigcup_{j \in \mathbb{Z}} S_j$. For this, it is enough to show that the action on $\mathbb{H} \times \mathbb{H}/M$ is free. Since
$$(A_r)^n \cdot (z_1, z_2) = (z_1, z_2) \Leftrightarrow ((A_r)^n z_1, (\bar{A}_r)^n z_2) = (z_1, z_2) \Leftrightarrow n = 0,$$
the action of $V$ on $\mathbb{H} \times \mathbb{H}/M$ is free. Hence the action is free on $Y - \bigcup_{j \in \mathbb{Z}} S_j$. Now we must show that the action is free on $\bigcup_{j \in \mathbb{Z}} S_j$. By $(A_r)^n$ a point $p$ on $S_k$ is mapped to a point on $S_{j+nr}$. If it is fixed, then $S_j \cap S_{j+nr} \neq \phi$ hence either $n = 1$ and $r = 1$ or $n = 0$. We use reduction to absurdity. Let us assume $n = 1$ and $r = 1$. By $A_1 \cdot p \in S_{j+1}$ and $A_1 \cdot p = p$, it follows that $p \in S_{j+1}$. Similarly, $p \in S_{j+2}$. Then $p \in S_j \cap S_{j+1} \cap S_{j+2}$. It is a contradiction to $S_j \cap S_{j+1} \cap S_{j+2} = \phi$. Therefore $n = 0$ and the action is free on $\bigcup_{j \in \mathbb{Z}} S_j$. To prove that $V$ is properly discontinuous, we show that for points $p, q \in Y^+$ there exist neighborhoods $U_1 \ni p, U_2 \ni q$ such that $gU_1 \cap U_2 \neq \phi$ for only finitely many $g \in V$. Since $V$ acts properly discontinuously on $\mathbb{H} \times \mathbb{H}/M$ and $\bigcup_{j \in \mathbb{Z}} S_j \subset Y^+$ is closed, this is clear if $p \notin \bigcup_{j \in \mathbb{Z}} S_j$ and $q \notin \bigcup_{j \in \mathbb{Z}} S_j$.

In the case where $p \in \bigcup_{j \in \mathbb{Z}} S_j$ and $q \notin \bigcup_{j \in \mathbb{Z}} S_j$, we can separate $p$ and $q$ by using the function $\rho(z_1, z_2) = \mathrm{Im} z_1 \cdot \mathrm{Im} z_2$ on $\mathbb{H} \times \mathbb{H}/M$. We put

$$U_2 = \{u | u \in Y^+, \rho \circ \Phi(u) < \rho \circ \Phi(q) + 1\}, \quad U_1 = Y^+ \backslash \overline{U}_2$$

then $p \in U_1, q \in U_2$ and $U_1 \cap U_2 = \phi$. Since $gU_1 = U_1 (\forall g \in V)$, $gU_1 \cap U_2 = \phi (\forall g \in V)$. In the case where $p, q \in \bigcup_{j \in \mathbb{Z}} S_j$, without loss of generality, we may assume $p \in S_k, q \in S_j$ with $k > j$. We put

$$\begin{aligned} U_1 &= \left\{ |u_k| < \varepsilon^{-\frac{1}{M}}, |v_k| < \varepsilon \right\}, \\ U_2 &= \left\{ |u_j| < \varepsilon^{-\frac{1}{M}}, |v_j| < \varepsilon \right\}, \end{aligned}$$

where $M = \max b_k = \max\{b_0, \cdots, b_{r-1}\}$. Take $\varepsilon$ sufficiently small, and $p \in U_1, q \in U_2$. Let us show that if $(A_r)^n U_1 \cap U_2 \neq \phi$, it is necessary that $j - k \leq n < 0$. First, we show that $U_1 \cap U_2 = \phi$. We use reduction to absurdity. Assume that $U_1 \cap U_2 \neq \phi$ and let $(u_k, v_k)$ and $(u_j, v_j)$ be the coordinates of $p \in U_1 \cap U_2$ in $R_k$ and $R_j$, respectively. Then

$$\begin{aligned} u_k &= u_j^a v_j^b, \\ v_k &= \frac{1}{u_j^c v_j^d}, \end{aligned}$$

where

$$\begin{pmatrix} a & b \\ -c & -d \end{pmatrix} = \begin{pmatrix} b_{k-1} & 1 \\ -1 & 0 \end{pmatrix} \cdots \begin{pmatrix} b_j & 1 \\ -1 & 0 \end{pmatrix}.$$

Hence if $k = j + 1$, $c = 1, d = 0$ and if $k > j + 1$

$$\frac{c}{d} = b_j - \cfrac{1}{b_{j+1} - \cfrac{1}{b_{j+2} - \cdots - \cfrac{1}{b_{k-2}}}} < b_j \leq M.$$

The second equation of coordinate transformation means $v_k u_j^c v_j^d = 1$, hence

$$|v_k| |u_j|^c |v_j|^d = 1.$$

On the other hand,

$$|v_k| < \varepsilon, |u_j| < \varepsilon^{-\frac{1}{M}}, |v_j| < \varepsilon$$

by $p \in U_1 \cap U_2$. Hence,

$$1 = |v_k| |u_j|^c |v_j|^d < \varepsilon \cdot \varepsilon^{-\frac{c}{M}} \cdot \varepsilon^d = \varepsilon^{1+d-\frac{c}{M}}, \quad \text{where } d - \frac{c}{M} = \frac{d}{M}(M - \frac{c}{d}) > 0.$$

This is not true for $\varepsilon < 1$. Therefore, $U_1 \cap U_2 = \phi$.
Similarly in the case $n > 0$, $(A_r)^n U_1 \cap U_2 = \phi$. Finally we consider the case $n < 0$ and $j > k + n$. Let us assume that $(A_r)^n U_1 \cap U_2 \neq \phi$. This time

$$(A_r)^n \cdot (u_k, v_k) = (u_k, v_k) \in R_{k+nr}, \quad j > k + nr$$

and thus

$$u_j = u_k^a v_k^b, \quad v_j = \frac{1}{u_k^c v_k^d},$$

where
$$\begin{pmatrix} a & b \\ -c & -d \end{pmatrix} = \begin{pmatrix} b_{j-1} & 1 \\ -1 & 0 \end{pmatrix} \cdots \begin{pmatrix} b_{k+nr} & 1 \\ -1 & 0 \end{pmatrix}.$$

Then estimating $|v_j||u_k|^c|v_k|^d$ leads to a contradiction, so $(A_r)^n U_1 \cap U_2 = \phi$! %
Hence if $(A_r)^n U_1 \cap U_2 \neq \phi$, it is necessary that $j - k \leq n < 0$ and there are only finitely many integers satisfying this condition. Therefore, the action $V$ on $Y^+$ is properly discontinuous. □

Hence, the two quotient spaces $Y^+/V$ and $(\mathbb{H} \times \mathbb{H}/M)/V$ are complex manifolds. We define the action of $G(M,V) = \left\{ \begin{pmatrix} \varepsilon & a \\ 0 & 1 \end{pmatrix} \mid a \in M, \varepsilon \in V \right\}$ on $\mathbb{H} \times \mathbb{H}$ such that

$$\begin{pmatrix} \varepsilon & a \\ 0 & 1 \end{pmatrix} \cdot (z_1, z_2) = (\varepsilon z_1 + a, \varepsilon^{-1} z_2 + \bar{a}).$$

Then $(\mathbb{H} \times \mathbb{H}/M)/V \cong \mathbb{H} \times \mathbb{H}/G(M,V)$. Indeed, for any $b \in M$ and $\varepsilon \in V$

$$\varepsilon \cdot (z_1 + b, z_2 + \bar{b}) = (\varepsilon z_1 + (\varepsilon b), \varepsilon^{-1} z_2 + (\overline{\varepsilon b})) = \begin{pmatrix} \varepsilon & \varepsilon b \\ 0 & 1 \end{pmatrix} \cdot (z_1, z_2).$$

From now on, we represent it $\mathbb{H} \times \mathbb{H}/G(M,V)$. In $Y^+$, the action of $A_r \in V$ sends a point $(u_k, v_k) \in S_k$ to a point $(u_k, v_k) \in S_{k+r}$. Therefore, $S_k$ corresponds to $S_{k+r}$ and non-singular rational curves $S_0, \cdots, S_{r-1}$ form a cycle in $Y(b_0, b_1, \cdots, b_{r-1}) = Y^+/V$. Hence it is expected that $Y(b_0, b_1, \cdots, b_{r-1})$ is a resolution of a cusp singularity. Indeed, the intersection matrix $(S_i \cdot S_j)$ is negative definite and satisfies Grauert's criterion.

**Theorem 3.11** (Grauert[2])**.** *Let $X$ be a non-singular complex surface. If the divisor $E$ in $X$ is such that the intersection matrix $A$ is negative definite, then we can blow down $E$ analytically; we get a normal complex surface $V$, in general with a singularity at the image 0 of $E$, and the projection $\pi : X \to V$ is a good resolution of $(V, 0)$ with exceptional divisor $E$.*

Hence $S_0, \cdots, S_{r-1}$ can be blown down to give an isolated normal point $p$ in a complex space $\overline{Y}(b_0, b_1, \cdots, b_{r-1})$. We have a holomorphic map

$$\pi : Y(b_0, b_1, \cdots, b_{r-1}) \to \overline{Y}(b_0, b_1, \cdots, b_{r-1}) \;,\; p = \pi(\bigcup_{j=0}^{r-1} S_j).$$

$\pi$ is a minimal resolution since the exceptional set contains no non-singular rational curves with self-intersection $-1$. Hence $p$ is a cusp singularity.
On the other hand, we consider $\overline{\mathbb{H} \times \mathbb{H}/G(M,V)}$ which is the completion of $\mathbb{H} \times \mathbb{H}/G(M,V)$ by adding the one point $\infty$. A function $f$ is said to be holomorphic at $\infty$ if there is a neighborhood $U$ of $\infty$ such that $f$ is holomorphic in $U \setminus \{\infty\}$ and continuous at $\infty$.

**Proposition 3.12.** $\overline{\mathbb{H} \times \mathbb{H}/G(M,V)}$ *is a normal complex space.*

*Proof.* It is enough to check Cartan's condition that there is some neighborhood $U$ of $\infty$ such that for any two different points $p_1, p_2 \in U - \{\infty\}$ there exists a holomorphic function $f$ in $U - \{\infty\}$ with $f(p_1) \neq f(p_2)$.
Since $\overline{Y}(b_0, b_1, \cdots, b_{r-1})$ is normal by Theorem 3.11, there is a neighborhood $V \ni p$ such that any two different points in $V$ can be separated by a holomorphic function in $V - \{p\}$. Let $\overline{\Phi} : \overline{Y}(b_0, b_1, \cdots, b_{r-1}) - \{p\} \to \mathbb{H} \times \mathbb{H}/G(M,V)$ be the map induced by the biholomorphic map $\Phi : Y - \bigcup_{j \in \mathbb{Z}} S_j \to \mathbb{C}^2/M$, and extend it over $\overline{Y}(b_0, b_1, \cdots, b_{r-1})$ with $\overline{\Phi}(p) = \infty$. Since $\overline{\Phi}$ is biholomorphic in $\overline{Y}(b_0, b_1, \cdots, b_{r-1}) - \{p\}$, $\overline{\Phi}$ transforms a holomorphic function in $V - \{p\}$ into a holomorphic function in $U - \{\infty\}$ with $U = \overline{\Phi}(V)$. Hence, $p_1, p_2$ can be separated by a holomorphic function in $U - \{\infty\}$. □

$\overline{\Phi} : \overline{Y}(b_0, b_1, \cdots, b_{r-1}) \to \overline{\mathbb{H} \times \mathbb{H}/G(M,V)}$ is a biholomorphic map, and $\infty$ is a cusp singularity. We put

$$W(d) = \{(z_1, z_2) \in \mathbb{H} \times \mathbb{H} \mid \mathrm{Im}(z_1) \cdot \mathrm{Im}(z_2) \geq d\},$$

then $\partial W(d)/G(M,V)$ is a link of $\infty$. The preimage of this link by the biholomorphic map above is also the link of the cusp singularity, $p \in \overline{Y}(b_0, b_1, \cdots, b_{r-1})$.
The definition of the link of a singularity is as follows.

**Definition 3.13.** *Let $X$ be a complex analytic variety of dimension at least $2$ with only one singular point $x$. Let $r : X \to [0, \infty)$ be a real analytic function on $X$ such that $r^{-1}(0) = x$, and $r$ is strictly pluri-subharmonic on $X - \{x\}$. For $\varepsilon > 0$ sufficiently small, $r$ has no critical points in $r^{-1}(0, \varepsilon]$. We define a link of germ $(X, x)$ to be $r^{-1}(\varepsilon)$. As we see in Lemma 3.17, any two links of $(X, x)$ are diffeomorphic.*

**Lemma 3.14.** *The function $\varphi(z_1, z_2) = \frac{1}{\mathrm{Im}(z_1) \cdot \mathrm{Im}(z_2)}$ on $\mathbb{H} \times \mathbb{H}$ is strictly pluri-subharmonic, and induces $\tilde{\varphi}$ which is strictly pluri-subharmonic on $\mathbb{H} \times \mathbb{H}/G(M, V)$.*

*Proof.* For the former part, we have to check that Levi matrix $L_\varphi$ is positive definite. Put $z_1 = x_1 + iy_1, z_2 = x_2 + iy_2$, then

$$\frac{\partial}{\partial z_j} = \frac{1}{2}(\frac{\partial}{\partial x_j} + \frac{1}{i}\frac{\partial}{\partial y_j}) \ , \ \frac{\partial}{\partial \bar{z}_j} = \frac{1}{2}(\frac{\partial}{\partial x_j} - \frac{1}{i}\frac{\partial}{\partial y_j}).$$

Since $\varphi = \frac{1}{y_1 y_2}$,

$$\frac{\partial^2 \varphi}{\partial z_1 \partial \bar{z}_1} = \frac{1}{2y_1^3 y_2} \ , \ \frac{\partial^2 \varphi}{\partial z_1 \partial \bar{z}_2} = \frac{1}{4y_1^2 y_2^2} \ ,$$

$$\frac{\partial^2 \varphi}{\partial z_2 \partial \bar{z}_1} = \frac{1}{4y_1^2 y_2^2} \ , \ \frac{\partial^2 \varphi}{\partial z_2 \partial \bar{z}_2} = \frac{1}{2y_1 y_2^3}.$$

Therefore,

$$4L_\varphi = \begin{pmatrix} \frac{2}{y_1^3 y_2} & \frac{1}{y_1^2 y_2^2} \\ \frac{1}{y_1^2 y_2^2} & \frac{2}{y_1 y_2^3} \end{pmatrix} = \frac{1}{y_1^3 y_2^3}\begin{pmatrix} 2y_2^2 & y_1 y_2 \\ y_1 y_2 & 2y_1^2 \end{pmatrix}.$$

It is easy to check that this matrix is positive definite.
For the latter part, we have to check that $\varphi$ is $G(M,V)$-invariant.
The action of $G(M,V) = \left\{ \begin{pmatrix} \varepsilon & a \\ 0 & 1 \end{pmatrix} \mid a \in M, \varepsilon \in V \right\}$ is

$$\begin{pmatrix} \varepsilon & a \\ 0 & 1 \end{pmatrix} \cdot (z_1, z_2) = (\varepsilon z_1 + a, \varepsilon^{-1} z_2 + \bar{a}),$$

and using the coordinates $(x_1, y_1, x_2, y_2)$, since $a, \bar{a} \in \mathbb{R}$,

$$\begin{pmatrix} \varepsilon & a \\ 0 & 1 \end{pmatrix} \cdot (x_1, y_1, x_2, y_2) = (\varepsilon x_1 + a, \varepsilon y_1, \varepsilon^{-1} x_2 + \bar{a}, \varepsilon^{-1} y_2).$$

Hence,
$$\begin{pmatrix} \varepsilon & a \\ 0 & 1 \end{pmatrix} \cdot \varphi = \frac{1}{(\varepsilon y_1)(\varepsilon^{-1} y_2)} = \frac{1}{y_1 y_2} = \varphi.$$

Indeed, $\varphi$ is $G(M,V)$-invariant and $\tilde{\varphi}$ is induced. $\square$

Moreover, define $\tilde{\varphi}(\infty) = 0$ and $\tilde{\varphi}$ is a real analytic function on $\overline{\mathbb{H} \times \mathbb{H}/G(M,V)}$. Hence, $\partial W(d)/G(M,V) = \tilde{\varphi}^{-1}(d)$ is the link of a cusp singularity $\infty$.

3.3. **Our results.** Let $J$ be the standard complex structure on $\mathbb{H} \times \mathbb{H}$. For $\varphi$ of Lemma 3.14, we put
$$\lambda = -J^* d\varphi \ , \ \omega = -dJ^* d\varphi \ , \ g(u,v) = \omega(u, Jv).$$

Then, $\omega$ is a symplectic form on $\mathbb{H} \times \mathbb{H}$ compatible with $J$, and $g$ is a $J$-invariant Riemannian metric. Moreover, $\alpha = \lambda \mid \partial W$ is a contact form on $\partial W$, and the Reeb vector field is $X_\alpha = J(Z)$ with $Z = \frac{\nabla \varphi}{\|\nabla \varphi\|^2}$. Since $\varphi$ is $G(M,V)$-invariant, these structures are $G(M,V)$-invariant. Therefore, $(\partial W/G, \tilde{\alpha})$ is a contact structure.

**Theorem 3.15.** $(\partial W/G, \tilde{\alpha})$ *is contactomorphic to a* Sol-*manifold* $(T_A, \beta_+ + \beta_-)$.

*Proof.* We put
$$A = \left\{ \begin{pmatrix} b_{r-1} & -1 \\ 1 & 0 \end{pmatrix} \begin{pmatrix} b_{r-2} & -1 \\ 1 & 0 \end{pmatrix} \cdots \begin{pmatrix} b_1 & -1 \\ 1 & 0 \end{pmatrix} \begin{pmatrix} b_0 & -1 \\ 1 & 0 \end{pmatrix} \right\}^{-1}.$$

Then, by the computation in the proof of Proposition 3.9,
$$A \begin{pmatrix} 1 \\ w_0 \end{pmatrix} = A_r^{-1} \begin{pmatrix} 1 \\ w_0 \end{pmatrix}, \ A \begin{pmatrix} 1 \\ \bar{w}_0 \end{pmatrix} = A_r \begin{pmatrix} 1 \\ \bar{w}_0 \end{pmatrix}, \text{ where } A_r < 1.$$

Hence we can put
$$\begin{pmatrix} 1 \\ w_0 \end{pmatrix} = \sqrt{w_0 - \bar{w}_0} v_+, \begin{pmatrix} 1 \\ \bar{w}_0 \end{pmatrix} = -\sqrt{w_0 - \bar{w}_0} v_- \text{ and } A_r^{-1} = a.$$

The map
$$F: T_A \to \partial W(1)/G \ ; \ ((x,y), z) \mapsto (y - w_0 x, \ a^z, \ y - \bar{w}_0 x, \ a^{-z})$$

is the diffeomorphism. Since
$$F(A \begin{pmatrix} 1 \\ w_0 \end{pmatrix}, 0) = (0, 1, a(w_0 - \bar{w}_0), 1) \sim F(\begin{pmatrix} 1 \\ w_0 \end{pmatrix}, 1) = (0, a, (w_0 - \bar{w}_0), a^{-1}),$$
$$F(A \begin{pmatrix} 1 \\ \bar{w}_0 \end{pmatrix}, 0) = (a^{-1}(\bar{w}_0 - w_0), 1, 0, 1) \sim F(\begin{pmatrix} 1 \\ \bar{w}_0 \end{pmatrix}, 1) = (\bar{w}_0 - w_0, a, 0, a^{-1}),$$

$F$ is well-defined and
$$H: \partial W(1)/G \to T_A \ ; \ (x_1, y_1, x_2, y_2) \mapsto ((\frac{x_2 - x_1}{w_0 - \bar{w}_0}, \frac{w_0 x_2 - \bar{w}_0 x_1}{w_0 - \bar{w}_0}), \log_a y_1)$$

is the inverse of $F$. Therefore, $\partial W/G$ is diffeomorphic to $T_A$.
Next we show that $F$ is the contactomorphism $(T_A, \beta_+ + \beta_-) \cong (\partial W(1)/G, \tilde{\alpha})$. For simplicity, let $d = 1$. The 1-form $\tilde{\alpha}$ is represented as follows.
$$\lambda = -J^* d\varphi = \frac{1}{y_1^2 y_2^2}(y_2 dx_1 + y_1 dx_2) \ , \ \alpha = \lambda \mid \partial W(1)$$

and $\alpha = y_2 dx_1 + y_1 dx_2$. This 1-form is $G$-invariant and induces $\tilde{\alpha}$ on $\partial W(1)/G$. Let us check $F^*\tilde{\alpha} = \sqrt{w_0 - \bar{w}_0}(\beta_+ + \beta_-)$.

$$\begin{aligned} F^*\tilde{\alpha} &= F^*(y_2 dx_1 + y_1 dx_2) = a^{-z} d(y - w_0 x) + a^z d(y - \bar{w}_0 x), \\ d(y - w_0 x) &= \sqrt{w_0 - \bar{w}_0} \cdot dx \wedge dy(v_+, \,\cdot\,), \\ d(y - \bar{w}_0 x) &= -\sqrt{w_0 - \bar{w}_0} \cdot dx \wedge dy(v_-, \,\cdot\,). \end{aligned}$$

Therefore,

$$\begin{aligned} F^*\tilde{\alpha} &= \sqrt{w_0 - \bar{w}_0}(a^{-z} dx \wedge dy(v_+, \,\cdot\,) - a^z dx \wedge dy(v_-, \,\cdot\,)) \\ &= \sqrt{w_0 - \bar{w}_0}(\beta_+ + \beta_-). \end{aligned}$$

$\square$

Thus, the induced contact structure on the link of a cusp singularity $\infty$ is contactomorphic to $(T_A, \beta_+ + \beta_-)$. This fact is also followed by the next remark.

**Remark 3.16.** *The Lie group $Sol^3$ acts freely on $\mathbb{H} \times \mathbb{H}$ as follows:*

$$(a, b; c) \cdot (z_1, z_2) = (e^c z_1 + a, e^{-c} z_2 + b).$$

*Thus we can identify $Sol^3$ with the orbit of a point, say $(i,i)$. It follows that the link $\partial W(1)/G(M,V)$ can be identified with $G(M,V)\backslash Sol^3$ and it is therefore a Sol-manifold. Moreover the complex structure $J$ and the strictly pluri-subharmonic function $\varphi$ are invariant under the action of $Sol^3$. Thus $\tilde{\alpha}$ is the contact form induced by the left invariant contact form $\alpha$ on $\partial W(1) \cong Sol^3$. By Remark 2.6, the contact structure is universally tight and uniquely determined. It is contactomorphic to $(T_A, \beta_+ + \beta_-)$, one of the pair consisting a bi-contact structure.*

The link of a singularity is unique in the following sense.

**Lemma 3.17.** *Let $(X, x)$ be a isolated surface singularity and $r_0$, $r_1$ be the different strictly pluri-subharmonic functions defining the links $L_0 = r_0^{-1}(\varepsilon)$, $L_1 = r_1^{-1}(\delta)$ of $(X, x)$. Then, restricting $\lambda_0 = -J^* dr_0$ and $\lambda_1 = -J^* dr_1$ to $L_0 = r_0^{-1}(\varepsilon)$ and $L_1 = r_1^{-1}(\delta)$ respectively, we get contact structures $(L_0, \alpha_0)$ and $(L_1, \alpha_1)$. For sufficient small $\varepsilon, \delta > 0$, these two contact structures are contactomorphic. In particular, the two links $L_0$ and $L_1$ are diffeomorphic.*

*Proof.* First, the contact structure $(L_0, \alpha_0)$ is independent of the choice of $\varepsilon$. Since $r_0$ is regular on $r_0^{-1}(0, \varepsilon]$, the gradient vector field $\nabla r_0$ does not vanish. Hence, the flow generated by the vector field $\frac{\nabla r_0}{\|\nabla r_0\|^2}$ connects any two links, and the isotopy induces a one parameter family of contact structures on $L_0 = r_0^{-1}(\varepsilon)$. Hence $(L_0, \alpha_0)$ is independent of the choice of $\varepsilon$ by the Gray stability. Now let us show the contactomorphism $(L_0, \alpha_0) \cong (L_1, \alpha_1)$. We put $N = r_0^{-1}((0, \varepsilon]) \cap r_1^{-1}((0, \delta])$ and take $\varepsilon_1, \delta_1 > 0$ such that $r_0^{-1}((0, \varepsilon_1]) \cup r_1^{-1}((0, \delta_1]) \subset N$. For any $\theta \in (0, \varepsilon_1] \cap (0, \delta_1]$, we define the link $L_0 = r_0^{-1}(\theta)$.

**Step 1.** $L_0 \cong L_1$(diffeomorphic).
We define $r_t = tr_1 + (1-t)r_0$. $r_t$ is a strictly pluri-subharmonic function and has no critical points on $N$. Hence, $\nabla r_t$ does not vanish on $N$. Put $v_t = (r_0 - r_1)\frac{\nabla r_t}{\|\nabla r_t\|^2}$. Integrating this time dependent vector field, and the image of $L_0$ by the isotopy $\rho_t$ is the link defined by $r_t$. Indeed,

$$\begin{aligned}
\frac{d}{dt}(\rho_t^* r_t) &= \rho_t^*(L_{v_t} r_t + \frac{dr_t}{dt}) \\
&= \rho_t^*(dr_t(v_t) + \frac{dr_t}{dt}) \\
&= \rho_t^*((r_0 - r_1) + (r_1 - r_0)) \\
&= 0
\end{aligned}$$

and $\rho_t$ maps $r_0^{-1}(\theta)$ to $r_t^{-1}(\theta)$. This means that $L_t = \rho_t(L_0) = r_t^{-1}(\theta)$ is the link defined by $r_t$. Since $L_0 \cong L_t$, in particular $L_0 \cong L_1$.

**Step 2.** $(L_0, \alpha_0) \cong (L_1, \alpha_1)$(contactomorphic).
The diffeomorphism $\rho_t : L_0 \to L_t$ pulls back $(L_t, \alpha_t)$ on $L_0$ to induce a contact structure $(L_0, \rho_t^* \alpha_t)$ for each $t \in [0, 1]$. Now we get a one parameter family of contact 1-forms on a closed manifold $L_0$. By the Gray stability, $(L_0, \alpha_0) \cong (L_0, \rho_t^* \alpha_t)$. Hence, $(L_0, \alpha_0) \cong (L_0, \rho_1^* \alpha_1) \cong (L_1, \alpha_1)$(contactomorphic). □

Since the singularity $(V, (0, 0, 0))$ in Theorem 3.5 is analytically equivalent to the cusp singularity $\infty$, the two links $(K, \xi_0|K)$ and $(\partial W/G, \tilde{\alpha})$ are contactomorphic by the above lemma.

**Theorem 3.18** (Main theorem). *$(K, \xi_0|K)$ is contactomorphic to $(T_A, \beta_+ + \beta_-)$.*

*Proof.* $(K, \xi_0|K)$ is equal to $(L_0, \alpha_0)$ of Lemma 3.17, where $(X, x) = (V, (0, 0, 0))$ and $r_0 = |x|^2 + |y|^2 + |z|^2$. Let $r_1$ be the strictly pluri-subharmonic function on $V$ induced by $\tilde{\varphi}$ of Lemma 3.14, then $(L_1, \alpha_1)$ is contactomorphic to $(T_A, \beta_+ + \beta_-)$ by Theorem 3.15. Since $(L_0, \alpha_0)$ and $(L_1, \alpha_1)$ are contactomorphic by Lemma 3.17, $(L_0, \alpha_0)$ is contactomorphic to $(T_A, \beta_+ + \beta_-)$. □

Next remark gives another proof of Theorem 3.18 by Honda's classification [4] and the theorem of Y.Lekili and B.Ozbagci [7].

**Remark 3.19.** *The canonical contact structure on the link of a cusp singularity is Stein fillable, in particular, tight and minimally twisting. Moreover, it is universally tight by the theorem of Y.Lekili and B.Ozbagci that Milnor fillable contact structures are universally tight [7]. On a Sol-manifold $T_A$, there is the unique positive contact structure that is universally tight and minimally twisting [4]. Therefore, $(K, \xi_0|K)$ is uniquely determined to be contactomorphic to $(T_A, \beta_+ + \beta_-)$.*

## 4. The relation to Mori's works

Atsuhide Mori proved the next theorem in [11]. Put

$$A_{m,k} = \begin{pmatrix} 1 & 0 \\ 1 & 1 \end{pmatrix} \begin{pmatrix} 1 & k_1 \\ 0 & 1 \end{pmatrix} \cdots \begin{pmatrix} 1 & 0 \\ 1 & 1 \end{pmatrix} \begin{pmatrix} 1 & k_m \\ 0 & 1 \end{pmatrix} \in SL(2; \mathbb{Z}),$$

where $m \in \mathbb{Z}_{>0}$, $k = (k_1, \cdots, k_m) \in (\mathbb{Z}_{\geq 0})^m$ and $k \neq 0$.

**Theorem 4.1** (Mori[11]). *The link of the isolated singular point $(0, 0, 0)$ of the hypersurface $f_{m,k}^{-1}(0) \subset \mathbb{C}^3$ is contactomorphic to $(T_{A_{m,k}}, \beta_+ + \beta_-)$, where*

$$f_{1,(k_1)}(x, y, z) = x^2 + (y - 2z^2)(y^2 + 2yz^2 + z^4 - z^{4+k_1}) \text{ and}$$
$$f_{2,(k_1,k_2)}(x, y, z) = x^2 + \{(y+z)^2 - z^{2+k_1}\}\{(y-z)^2 + z^{2+k_2}\}.$$

The polynomials $f_{1,(k_1)}$ and $f_{2,(k_1,k_2)}$ are corresponding to $x^2+y^3+z^{6+k_1}+xyz$ and $x^2+y^{4+k_1}+z^{4+k_2}+xyz$ of Theorem 1.1, respectively. The finite sequence $(k_1+2,\cdots,k_m+2) \in \mathbb{Z}_{\geq 2}^m$ is equal to the sequence $\mathbf{b}^*$ in the second section of [1]. Then Theorem 1.1 says that

$$f_{3,(k_1,k_2,k_3)}(x,y,z) = x^{3+k_1}+y^{3+k_2}+z^{3+k_3}+xyz.$$

It is the answer for Problem 4.6 of [11] for the case where $m = 3$. Moreover there is no surface singularity for the case where $m > 3$ because of Laufer and Karras's criterion (see Theorem 3.5).

**Remark 4.2.** If $k = 0$, though this case is excluded by the above condition $k \neq 0$, $T_{A_{m,0}}$ is $S^1$-bundle over $T^2$ with the Euler characteristic $-m$. It is called a Nil-manifold. The polynomials $f_{1,0} = x^2+y^3+z^6+xyz$, $f_{2,0} = x^2+y^4+z^4+xyz$, $f_{3,0} = x^3+y^3+z^3+xyz$ are corresponding to $x^2+y^3+z^6$, $x^2+y^4+z^4$, $x^3+y^3+z^3$. The associated isolated surface singularities are called $\tilde{E}_8, \tilde{E}_7, \tilde{E}_6$ (simple elliptic singularities) and their links are Brieskorn Nil-manifolds $M(2,3,6), M(2,4,4), M(3,3,3)$, respectively.

Moreover Mori also computed the Euler characteristic of the Milnor fiber $P_{m,k}$ associated with $f_{m,k}$ (Theorem 4.5 of [11]). The result

$$\chi(P_{1,(k_1)}) = 11+k_1,\ \chi(P_{2,(k_1,k_2)}) = 10+k_1+k_2,\ \chi(P_{3,(k_1,k_2,k_3)}) = 9+\sum_{i=1}^{3}k_i$$

is also seen in the second section of [1] and [5]. He got these results in the context of constructing the examples of Lutz-Mori twists, the generalization of Lutz twists to higher dimensions.

4.1. **Lutz-Mori twists.** This section is based on [11]. See also the last section of [8].

**Definition 4.3.** A Geiges pair on $M^{2n+1}$ is a pair of contact forms $\alpha_+$ and $\alpha_-$ such that for all $0 \leq k \leq n-1$

$$\alpha_+ \wedge d\alpha_+^k \wedge d\alpha_-^{n-k} = \alpha_- \wedge d\alpha_-^k \wedge d\alpha_+^{n-k} = 0,\ \ \alpha_+ \wedge d\alpha_+^n = -\alpha_- \wedge d\alpha_-^n > 0.$$

**Definition 4.4.** Suppose that $(N^{2n+1}, \alpha)$ is a contact manifold containing a closed codimension $2$ contact submanifold $(M^{2n-1}, \alpha_+)$ with trivial normal bundle such that $\alpha_+ = \alpha|M^{2n-1}$ is one of some Geiges pair $(\alpha_+, \alpha_-)$ on $M^{2n-1}$. We can take the polar coordinates $(r,\theta)$ such that $\ker\alpha = \ker(\alpha_+ + r^2 d\theta)$ on $\{r < \varepsilon\}$. Then the 1-form

$$\lambda = \frac{1-\cos s}{2}\alpha_+ + \frac{1+\cos s}{2}\alpha_- - s\sin s\,d\theta$$

defines the contact structure on $M^{2n-1} \times D^2$, where $(s,\theta)$ are the polar coordinates of the disk $D^2$ with the radius $\pi$. We call $(M^{2n-1} \times D^2, \lambda)$ the Lutz-Mori tube. Putting $s = r+\pi$ and slightly deform the contact structure $\ker(\alpha_+ + r^2 d\theta)$ on $\{0 < r < \varepsilon\}$, we can insert the Lutz-Mori tube along the contact submanifold $M^{2n-1}$ to get a possibly new contact structure on $N^{2n+1}$. We call this operation the Lutz-Mori twist.

**Example 4.5.** In the case where $n = 1$, each connected component of $M^1$ is the circle $S^1$ with a coordinate $z$. Then $(\alpha_+, \alpha_-) = (dz, -dz)$ is a Geiges pair on $S^1$.

*The contact form* $\lambda = -\cos s dz - s \sin s d\theta$ *defines the usual Lutz tube structure on the solid torus* $S^1 \times D^2$.

**Example 4.6.** *Suppose that* $M^3 = T_A$ $(n = 2)$. $TM^3$ *admits a frame* $(e_1, e_2, e_3)$ *with*

$$[e_3, e_2] = e_1, \ \ [e_3, e_1] = e_2, \ \ and \ \ [e_1, e_2] = 0,$$

*because* $M^3$ *is a Sol-manifold. Indeed,*

$$e_1 = \frac{a^z v_- + a^{-z} v_+}{2}, \ \ e_2 = \frac{a^z v_- - a^{-z} v_+}{2} \ \ and \ \ e_3 = \frac{1}{\log a}\frac{\partial}{\partial z}$$

*satisfy the conditions. The dual coframe* $(\alpha_1, \alpha_2, \alpha_3) = (\beta_+ + \beta_-, \beta_+ - \beta_-, \log a \cdot dz)$ *satisfies*

$$d\alpha_1 = \alpha_2 \wedge \alpha_3, \ \ d\alpha_2 = \alpha_1 \wedge \alpha_3 \ \ and \ \ d\alpha_3 = 0.$$

*Then* $(\alpha_+, \alpha_-) = (\beta_+ + \beta_-, \beta_+ - \beta_-)$ *is a Geiges pair on* $T_A$. *The contact form*

$$\lambda = \beta_+ - \cos s \cdot \beta_- - s \sin s \cdot d\theta$$

*defines the Lutz-Mori tube* $(T_A \times D^2, \lambda)$.

Performing the Lutz-Mori twist, the resultant contact structure $\xi$ on $N^{2n+1}$ contains a plastikstufe, the generalization of an overtwisted disk to higher dimensions. Therefore $(N^{2n+1}, \xi)$ is not strongly fillable (see [12]). In the case where $(N^5, \ker \alpha) = (S^5_\varepsilon, \xi_0)$ and $(M^3, \ker \alpha_+) = (K, \xi_0|K)$ of Theorem 3.5, the resultant contact structure $(S^5_\varepsilon, \xi)$ is an exotic contact structure because it cannot be the boundary of the standard symplectic ball $(B^6, \omega_0)$. Moreover Mori proved that $\xi$ can be deformed via contact structures to a spinnable foliation (Theorem 3.7 of [11]).

## References

[1] W.Ebeling, C.T.C.Wall. *Kodaira singularities and an extension of Arnold's strange duality.* Compositio Math.,56:3-77,1985.
[2] H.Grauert. *Über Modifikationen und exzeptionelle analytische Mengen.* Math.Ann.,146:331-368, 1962.
[3] F.Hirzebruch. *Hilbert modular surfaces.* L'Enseignement Math.,19:183-281,1973.
[4] K.Honda. *On the classification of tight contact structures II.* J. Diff. Geom., 55:83-143, 2000.
[5] U.Karras. *Deformations of cusp singularities.* Proc. Symp. Pure Math. XXX, ed. R.O. Wells jr. pp. 37-40. Providence, R.I.:Am. Math. Soc., 1977.
[6] H.B.Laufer. *Minimally elliptic singularities.* Am.J.Math.,99:1257-1295,1977.
[7] Y.Lekili, B.Ozbagci. *Milnor fillable contact structures are universally tight.* Math. Res. Lett. 17 (2010), no. 6, 1055-1063.
[8] P.Massot, K.Niederkrüger, C.Wendl. *Weak and strong fillability of higher dimensional contact manifolds*, preprint(2011), arXiv:1111.6008v1[math.SG].
[9] J.Milnor. *Singular Points of Complex Hypersurfaces.* Annals of Maths. Study 61,Princeton University Press, Princeton, 1968.
[10] Y.Mitsumatsu. *Anosov flows and non-Stein symplectic manifolds.* Ann. Inst. Fourier,45:1407-1421,1995.
[11] A.Mori. *Reeb foliations on* $S^5$ *and contact 5-manifolds violating the Thuston-Bennequin inequality*, preprint(2009), arXiv:0906.3237v2[math.GT].
[12] K.Niederkrüger. *The plastikstufe - a generalization of overtwisted disk to higher dimensions.* Algebr.Geom.Topol.6:2473-2508, 2006.

Graduate School of Mathematical Sciences, University of Tokyo, 3-8-1 Komaba, Meguro-ku, Tokyo 153-8914, Japan.
*E-mail address*: nkasuya@ms.u-tokyo.ac.jp